\documentclass[sagev,times]{sagejny}

\usepackage{dcolumn}
\usepackage{fix-cm}
\usepackage{epstopdf}

\usepackage[hidelinks]{hyperref}
\usepackage{graphicx}
\usepackage{dcolumn}
\usepackage{subcaption} 
\usepackage{longtable}
\usepackage{arydshln}
\theoremstyle{plain}
\newtheorem{assumption}{Assumption}
\theoremstyle{plain}
\newtheorem{theorem}{Theorem}
\newcommand\independent{\protect\mathpalette{\protect\independenT}{\perp}}
\def\independenT#1#2{\mathrel{\rlap{$#1#2$}\mkern2mu{#1#2}}} 

\usepackage[a4paper]{geometry}

\makeatletter
\renewcommand\subsubsection{\@startsection{subsubsection}{3}{\z@}%
	{0.5\@bls plus .3\@bls minus .1\@bls}%
	{0.5em\@afterindentfalse}%
	{\sagesf\normalsize\itshape}}
\makeatother

\begin{document}
\runninghead{A. Lindmark et al.}
\title{Sensitivity analysis for unobserved confounding of direct and indirect effects using uncertainty intervals}
\author{Anita Lindmark\affilnum{a},
	Xavier de Luna\affilnum{a} and
	Marie Eriksson\affilnum{a}}

\affiliation{\affilnum{1}Department of Statistics, Ume\r{a} School of Business and Economics, Ume\r{a} University, Ume\r{a}, Sweden}
\corrauth{Anita Lindmark, Department of Statistics, Ume\r{a} School of Business and Economics, Ume\r{a} University, SE-901 87 Ume\r{a}, Sweden}

\email{anita.lindmark@umu.se}


\begin{abstract}
To estimate direct and indirect effects of an exposure on an outcome from observed data strong assumptions about unconfoundedness are required. Since these assumptions cannot be tested using the observed data, a mediation analysis should always be accompanied by a sensitivity analysis of the resulting estimates. In this article we propose a sensitivity analysis method for parametric estimation of direct and indirect effects when the exposure, mediator and outcome are all binary. The sensitivity parameters consist of the correlation between the error terms of the mediator and outcome models, the correlation between the error terms of the mediator model and the model for the exposure assignment mechanism, and the correlation between the error terms of the exposure assignment and outcome models. These correlations are incorporated into the estimation of the model parameters and identification sets are then obtained for the direct and indirect effects for a range of plausible correlation values. We take the sampling variability into account through the construction of uncertainty intervals. The proposed method is able to assess sensitivity to both mediator-outcome confounding and confounding involving the exposure. To illustrate the method we apply it to a mediation study based on data from the Swedish Stroke Register (Riksstroke). 
\end{abstract}
\maketitle

\section{Introduction}
Evaluating the effect of an exposure (or treatment) on an outcome is a common goal in medical studies, e.g. clinical trials, epidemiological studies, and quality of care evaluations. Mediation analysis seeks to further decompose this effect into direct and indirect effects (i.e. effects that take the pathway through some intermediate variable, a mediator), in order to better understand the causal mechanisms at work and where best to target interventions. As an example, there is evidence that patients who live alone tend to have worse prognosis after stroke than those who cohabit \cite{lindmark2014,eriksson2009}. It is of interest to uncover the causal mechanisms behind this association. Is there some intermediate variable affected by cohabitation status that in turn has an effect on the outcome after stroke?  Is this effect to a large extent due to structural problems within the care system or could there be some property of the patients themselves that leads to a less advantageous outcome?

Traditional approaches for estimating mediation effects have relied on parametric linear regression models, through e.g. the ``product method" popularized by Baron and Kenny \cite{baron1986}. Formalizing the concepts of direct and indirect effects in the causal inference framework has led to methodological developments and better understanding of which assumptions are required for estimation of and inference about these effects \cite{pearl2014}. The traditional parametric methods have been generalized to allow for exposure-mediator interactions and a broader class of mediator and outcome types \cite{vanderweele2009,valeri2013,vanderweele2010}. Non-parametric and semi-parametric estimation methods have also been proposed \cite{tchetgen2012,huber2014,tenhave2012,daniels2012,vanderweele2009b,goetgeluk2008}. 

To estimate direct and indirect effects from observed data strong assumptions are made about unconfoundedness of the relationships between exposure, mediator, and outcome, assumptions that are not testable using the observed data. In theory one could safeguard against confounding involving the exposure by randomizing it, but even in situations where this is possible it is difficult to rule out confounding between the mediator and the outcome. In observational studies, where randomization is not possible, a solution is to adjust for observed pre-exposure covariates, i.e. covariates that either temporally precede the exposure or through e.g. subject-matter knowledge are guaranteed to be unaffected by it. In situations where the exposure affects confounders of the mediator-outcome relation, additional assumptions are required \cite{robins1992,avin2005,petersen2006,destavola2015}. In any case, because unobserved confounding can seldom be ruled out, a sensitivity analysis of the effect of its existence is an essential complement to a mediation analysis adjusted for observed confounders. 

The sensitivity analysis techniques that have been suggested in the literature have up to our knowledge  exclusively focused on mediator-outcome confounding \cite{vanderweele2010S,imai2010b,tchetgen2012,daniels2012}. In the parametric modeling setting of this paper, VanderWeele \cite{vanderweele2010S} suggests specifying a bias factor which is then used to correct estimates and confidence intervals. This requires specification of the effect of the unknown confounder on the outcome given exposure, mediator, and observed covariates as well as the relation between the exposure and the unmeasured confounder given the mediator and observed covariates. This may be difficult to do in practice; see also Hafeman \cite{hafeman2011}. Recently, le Cessie \cite{lecessie2016} proposed specifying the effect of the unobserved confounder on the (continuous) mediator and on the (continuous, binary or count) outcome  directly in the parametric regression models. Bias formulas were derived under the assumption that the unobserved confounder follows a normal distribution.  Imai et al. \cite{imai2010b} have proposed an alternative method, also model based, that uses the fact that unobserved confounding of the mediator-outcome relation induces a correlation between the error terms in the parametric regression models for the mediator and outcome. They then derive expressions for the direct and indirect effects that takes this correlation into account. This has the advantage that only one sensitivity parameter needs to be specified. This method has been implemented through the function \texttt{medsens} in the R \cite{r2015} package \textbf{\texttt{mediation}} \cite{tingley2013,tingley2014}, for continuous mediators and outcomes and for situations where either the  mediator or the outcome is binary. The models available are linear regression for continuous variables and binary probit regression. However, in the current implementation, if a binary outcome model is used it cannot include an exposure-mediator interaction term, which is often important in order to fully capture the dynamics of mediation \cite{vanderweele2015}. 

In this article we propose a sensitivity analysis that allows us to investigate unobserved mediator-outcome confounding \emph{and} unobserved confounding of the exposure-mediator and exposure-outcome relations, when parametric models are used to obtain estimators for conditional and marginal direct and indirect effects. Building on a proposal by Genb\"{a}ck et al. \cite{genback2014} for sensitivity analysis of regression parameters in the presence of non-ignorable missing outcomes, the sensitivity analysis introduced here is based on correlations between the error terms of the exposure, mediator, and outcome models. These correlations are incorporated in the estimation of the regression parameters, upon which the direct and indirect effects estimates are based, through a modified maximum likelihood (ML) procedure. Sampling variability is further taken into account through the construction of uncertainty intervals (UIs) \cite{vansteelandt2006}. We present the sensitivity analysis for binary mediator and outcome variables, although the same ideas can be used for continuous mediators and outcomes. Finally, our approach allows for an exposure-mediator interaction term in the outcome model, in contrast with some of the existing methods described above. To illustrate the use of the sensitivity analysis proposed, it was applied to a study using data from the Swedish national quality register for stroke. We investigated the effect of living alone on the probability of death or dependency in activities of daily living (ADL) 3 months after stroke among male patients registered in 2012, and to which extent this effect was mediated by stroke severity, the theory being that patients living alone are less likely to recognize stroke symptoms and therefore arrive to the hospital later and with a more severe stroke than patients who cohabit. We then used the proposed sensitivity analysis technique to assess the sensitivity of our findings to unobserved confounding.

This article is organized as follows. In Section \ref{sec:causaleffects} we first introduce mediation analysis using  counterfactuals, followed by a definition of direct and indirect causal effects. Section \ref{sec:identification} presents assumptions necessary to identify these effects from observed data and gives the main identification result. A parametric estimation approach using regression models is described and we suggest probit based estimators for binary outcomes and mediators. In Section \ref{sec:sensanalys} we introduce a new method for sensitivity analysis to unobserved confounding and in Section \ref{sec:casestudy} the latter is applied to a mediation study. Finally, the study and results are summarized in Section \ref{sec:discussion}.

\section{Causal effects in mediation analysis}
\label{sec:causaleffects}
\begin{figure}[ht]
	\centering	
	\includegraphics[]{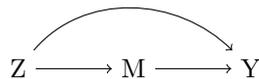}
	\caption{A directed acyclic graph showing the relationships between exposure $Z$, mediator $M$, and outcome $Y$.}
	\label{fig:DAG}
\end{figure}

Let $Z_i$ be an exposure such that $Z_i=1$ if individual $i$ is exposed, and $0$ otherwise. Let $Y_i$ be an outcome variable, and suppose further that we have an intermediate variable, $M_i$, on the pathway between the exposure and the outcome (see the directed acyclic graph \cite{lauritzen1996} in Figure~\ref{fig:DAG}). We refer to $M_i$ as a \emph{mediator} of the relationship between $Z_i$ and $Y_i$. In mediation analysis the goal is to distinguish effects that are \emph{indirect}, i.e. where the exposure affects some intermediate variable (or variables) of interest, and this intermediate variable in turn affects the outcome, from effects that are \emph{direct}, i.e. the effect of the exposure on the outcome not mediated through this intermediate variable. In our setting the indirect effect corresponds to  the path from $Z$ to $Y$ that passes through $M$, and the direct effect corresponds to the path from $Z$ to $Y$ that does not pass through $M$. To define these causal effects formally we will use counterfactuals, as formulated by Robins and Greenland \cite{robins1992} and Pearl \cite{pearl2001} for mediation, see also VanderWeele \cite{vanderweele2015}. 

\subsection{Counterfactuals}
Let us start by ignoring the role of $M$ and instead focus on the so called \emph{total effect} of $Z$ on $Y$. For each individual $i$ we would like to contrast the outcome had this individual been exposed to the outcome if the individual had not been exposed. To this end we introduce $Y_i(z)$, the potential outcome for individual $i$ under exposure level $z$.  The desired contrast would then be given by $Y_i(1)-Y_i(0)$. However, since only one of these outcomes can be observed for each individual we often seek to estimate the average causal effect of $Z$ on $Y$ over the entire population, giving the definition of the (average) total effect as:
$$TE= \mathbb{E}\left[Y_i\left(1\right)-Y_i\left(0\right)\right].$$

Returning to mediation analysis, we need to expand these counterfactuals to take into account the role of $M$. We let $M_i(z)$ denote the potential value of the mediator for individual $i$ under exposure level $z$. Further, since $Y$ is a function of both $Z$ and $M$, we denote the potential outcome under exposure level $z$ and mediator level $m$ as $Y_i(z,m)$. In addition we can express the composite potential outcome if the exposure $Z_i$ were set to the value $z$ and the mediator $M_i$ were set to its value under exposure level $Z_i=z'$: $Y(z,M(z'))$.

\subsection{Definition of direct and indirect effects}
\label{sec:definitions}
There are different definitions of direct and indirect effects \cite{robins1992,vanderweele2015}. Here we focus on the two most commonly defined effects, expressed on the difference scale. The \emph{natural direct effect}, $NDE$, sometimes referred to as the ``pure direct effect" \cite{robins1992}, is defined as the effect of $Z$ on $Y$ when allowing the mediator to vary as it would naturally if all individuals in the population were unexposed 
$$NDE =\mathbb{E}\left[Y_i\left(1,M_i(0)\right)-Y_i\left(0,M_i(0)\right)\right].$$

The \emph{natural indirect effect}, $NIE$, sometimes referred to as the ``total indirect effect" \cite{robins1992}, is defined as the effect on $Y$ of changing the mediator from its potential value when $Z=1$, $M_i(1)$, to its potential value when $Z=0$, $M_i(0)$, while keeping $Z$ fixed at $Z=1$
$$NIE=\mathbb{E}\left[Y_i\left(1,M_i(1)\right)-Y_i\left(1,M_i(0)\right)\right].$$

The natural direct and indirect effects are of interest when describing and evaluating the causal mechanisms at work \cite{pearl2001}. An important property of the definition of the $NDE$ and $NIE$ using these counterfactual-based definitions is that the total effect decomposes into the sum of the $NDE$ and $NIE$, i.e.  $TE=NDE+NIE$ \cite{pearl2014}. 

Note that the total effect can be decomposed in different ways, leading to alternative definitions of the natural direct and indirect effect. The ``total direct effect" and ``pure indirect effect" are defined as $\mathbb{E}\left[Y_i\left(1,M_i(1)\right)-Y_i\left(0,M_i(1)\right)\right]$ and  $\mathbb{E}\left[Y_i\left(0,M_i(1)\right)-Y_i\left(0,M_i(0)\right)\right]$. 

\section{Identification and estimation of direct and indirect effects}
\label{sec:identification}
To identify direct and indirect effects from observed data, certain assumptions need to be fulfilled. First we need to make a consistency assumption which states that (i) for an individual $i$ with observed exposure $Z_i=z$ and observed mediator $M_i=m$ the observed outcome is given by $Y_i=Y_i(z,m)$, (ii) for an individual $i$ with observed exposure $Z_i=z$ the observed mediator is given by $M_i=M_i(z)$, and (iii) for an individual $i$ with observed exposure $Z_i=z$ the observed outcome is given by $Y_i=Y_i(z,M(z))$ \cite{vanderweele2009,vansteelandt2012}. We also need to make a no interference assumption, meaning that the exposure level of one individual does not have an effect on the mediator or the outcome of another individual \cite{destavola2015}.

In addition to consistency and no interference we need to make assumptions about confounding. There are different formulations of these assumptions \cite{pearl2014}, here we use the \emph{sequential ignorability} assumption formulated by Imai et al.\cite{imai2010b,imai2010} 

\begin{assumption}
	\label{ass:seqign}
	\emph{Sequential ignorability (Imai et al. \cite{imai2010b}).} There exists a set of observed covariates $\boldsymbol{X}$ such that
	\begin{align}
	&\left\lbrace Y_i(z',m),M_i(z)\right\rbrace  \independent Z_i|\boldsymbol{X}_i=\boldsymbol{x}, \label{eg:seqign1}\\
	&Y_i(z',m)\independent M_i(z)|Z_i=z,\boldsymbol{X}_i=\boldsymbol{x}, \label{eg:seqign2}
	\end{align}
	where $0<P\left(Z_i=z|\boldsymbol{X}_i=\boldsymbol{x} \right) $ and $0<P\left(M_i(z)=m|Z_i=z,\boldsymbol{X}_i=\boldsymbol{x} \right) $ for $z,z'=0,1$, and all $\boldsymbol{x}\in\mathcal{X}$ (where $\mathcal{X}$ is the support of $\boldsymbol{X}_i$) and all $m\in\mathcal{M}$ (where $\mathcal{M}$ is the support of $M$).
\end{assumption}

\begin{figure}[ht]
	\centering	
	\includegraphics[]{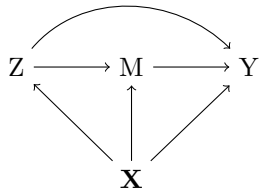}
	\caption{A directed acyclic graph showing the relationships between exposure $Z$, mediator $M$, outcome $Y$, and the set of observed confounders $\boldsymbol{X}$.}
	\label{fig:DAG2}
\end{figure}

Note that we use upper case letters do denote random variables and lower case letters to denote their realizations. The first part of this assumption says that there is no unobserved confounding of the exposure-mediator and exposure-outcome relationship given the observed covariates $\boldsymbol{X}_i$. The second part says that given $\boldsymbol{X}_i$ and the observed exposure $Z_i$ there is no confounding of the mediator-outcome relationship. Note that $\boldsymbol{X}_i$ need to be pre-exposure (i.e. not affected by the exposure) covariates, otherwise additional assumptions are required to identify the natural direct and indirect effects \cite{robins1992,avin2005,petersen2006,destavola2015}. Interpreting the DAG in Figure \ref{fig:DAG2} as a non-parametric structural equation model with independent error terms \cite{pearl2001}, it illustrates a situation where Assumption \ref{ass:seqign} is fulfilled. 

If consistency, no interference, and Assumption \ref{ass:seqign} are fulfilled the direct and indirect effects are identified through the following result\cite{pearl2001,pearl2009,imai2010b,vanderweele2009} 

\begin{theorem}
	\label{thm:medforms}
	(Pearl. \cite{pearl2001}) If Assumption \ref{ass:seqign} holds the natural direct and indirect effects conditional on the covariates $\boldsymbol{x}$ are identified as
	\begin{align}
	NDE\left( \boldsymbol{x}\right)  &=\sum_{m}\left[\mathbb{E}\left(Y_i|Z_i=1,M_i=m,\boldsymbol{X}_i=\boldsymbol{x} \right)- \mathbb{E}\left(Y_i|Z_i=0,M_i=m,\boldsymbol{X}_i=\boldsymbol{x} \right) \right]\times \notag\\
	&\qquad\quad P\left( M_i=m|Z_i=0,\boldsymbol{X}_i=\boldsymbol{x}\right) \label{eq:thm1nde},\\[7pt]
	NIE\left( \boldsymbol{x}\right) &=\sum_{m}\mathbb{E}\left(Y_i|Z_i=1,M_i=m,\boldsymbol{X}_i=\boldsymbol{x} \right)\times\notag\\
	&\qquad\quad\left[P\left( M_i=m|Z_i=1,\boldsymbol{X}_i=\boldsymbol{x}\right)- P\left( M_i=m|Z_i=0,\boldsymbol{X}_i=\boldsymbol{x}\right) \right]. \label{eq:thm1nie}
	\end{align}
\end{theorem}
For continuous mediators the sums in Theorem \ref{thm:medforms} are replaced by integrals and probabilities are replaced by densities. The $NDE$ and $NIE$ for the population (marginal effects) can be obtained by summing (or integrating) \eqref{eq:thm1nde} and \eqref{eq:thm1nie} over $\boldsymbol{x}$, e.g. the marginal natural indirect effect is given by $NDE =\sum_{m}\sum_{\boldsymbol{x}}\left[\mathbb{E}\left(Y_i|Z_i=1,M_i=m,\boldsymbol{X}_i=\boldsymbol{x} \right)- \mathbb{E}\left(Y_i|Z_i=0,M_i=m,\boldsymbol{X}_i=\boldsymbol{x} \right) \right]\times P\left( M_i=m|Z_i=0,\boldsymbol{X}_i=\boldsymbol{x}\right)P\left( \boldsymbol{X}_i=\boldsymbol{x}\right)$.  

Note that the corresponding identification result for the alternative definitions of the natural direct and indirect effects introduced at the end of Section \ref{sec:definitions} is obtained by replacing $P\left( M_i=m|Z_i=0,\boldsymbol{X}_i=\boldsymbol{x}\right)$ in \eqref{eq:thm1nde} with $P\left( M_i=m|Z_i=1,\boldsymbol{X}_i=\boldsymbol{x}\right)$ and $\mathbb{E}\left(Y_i|Z_i=1,M_i=m,\boldsymbol{X}_i=\boldsymbol{x} \right) $ in \eqref{eq:thm1nie} with $ \mathbb{E}\left(Y_i|Z_i=0,M_i=m,\boldsymbol{X}_i=\boldsymbol{x} \right) $.

\subsection{Parametric modeling and estimation}
The expressions in Theorem \ref{thm:medforms} are estimable from the observed data, both with and without specifying parametric models for the outcome and mediator. For a review of different estimation methods for direct and indirect effects, see Ten Have et al. \cite{tenhave2012}. Here the focus will be on sensitivity analysis for parametric estimation.

The classic  ``product method" approach \cite{judd1981} made popular by Baron and Kenny \cite{baron1986}, operates in the Linear Structural Equation Models (LSEM) framework, but predates the definition of mediation effects using counterfactual notation. It has been extended to allow for exposure-mediator interactions and binary mediators and outcomes \cite{valeri2013,vanderweele2009,vanderweele2010}. This approach estimates the natural direct and indirect effects by specifying parametric regression models for the outcome and the mediator. In the case of a continuous mediator and continuous outcome, the following linear regression models are fitted:
\begin{align}
&\mathbb{E}\left( M_i|Z_i=z,\boldsymbol{X}_i=\boldsymbol{x}\right) = \sum_{m}mP\left( M_i=m|z,\boldsymbol{x}\right) = \beta_0 + \beta_1z + \boldsymbol{\beta}_2^\top\boldsymbol{x} + \boldsymbol{\beta}_3^\top z\boldsymbol{x}, \label{eq:contM} \\
&\mathbb{E}\left( Y_i|Z_i=z,M_i=m,\boldsymbol{X}_i=\boldsymbol{x}\right) =\theta_0+\theta_1z + \theta_2m + \theta_3zm + \boldsymbol{\theta}_4^\top\boldsymbol{x}+\boldsymbol{\theta}_5^\top z\boldsymbol{x} + \boldsymbol{\theta}_6^\top m\boldsymbol{x} + \boldsymbol{\theta}_7^\top zm\boldsymbol{x}. \label{eq:contY}
\end{align}  
These models are often simplified to only include the main effects of the covariates $\boldsymbol{X}_i$. Substituting these simplified versions of the models into \eqref{eq:thm1nde}-\eqref{eq:thm1nie} yields expressions for the mediation effects in terms of the regression coefficients, $NDE(\boldsymbol{x})=\theta_1+\theta_3(\beta_0+\boldsymbol{\beta}_2^\top\boldsymbol{x})$ and $NIE(\boldsymbol{x})=\beta_1(\theta_2+\theta_3)$ \cite{vanderweele2015}. Corresponding expressions based on the more general models in \eqref{eq:contM} and \eqref{eq:contY} can be similarly derived. Given that the assumptions in the previous section are fulfilled and the regression models are correctly specified, these expressions yield consistent estimators of the mediation effects \cite{vanderweele2009,vanderweele2010}. 

When the outcome and mediator are both binary, expressions for the natural direct and indirect effects have been derived on the odds ratio scale using logistic regression models for the outcome and mediator \cite{valeri2013,vanderweele2010}. Here we develop alternative expressions based on probit regression models. Let us assume that $M_i$ and $Y_i$ are both binary random variables and can be modeled by $M_i = I(M_i^*>0)$, where
\begin{equation}
M_i^*=\beta_0 + \beta_1Z_i + \boldsymbol{\beta}_2^\top\boldsymbol{X}_i + \boldsymbol{\beta}_3^\top Z_i\boldsymbol{X}_i + \eta_i \label{mstar},
\end{equation}
and $Y_i = I(Y_i^*>0)$, where
\begin{equation}
Y_i^*=\theta_0+\theta_1Z_i + \theta_2M_i + \theta_3Z_iM_i + \boldsymbol{\theta}_4^\top\boldsymbol{X}_i + \boldsymbol{\theta}_5^\top Z_i\boldsymbol{X}_i + \boldsymbol{\theta}_6^\top M_i\boldsymbol{X}_i + \boldsymbol{\theta}_7^\top Z_iM_i\boldsymbol{X}_i + \xi_i \label{ystar}.
\end{equation}
$I\left( A>0\right)$ is an indicator variable that takes the value 1 if $A>0$ and 0 otherwise. The error terms $\eta_i$ and $\xi_i$ are both assumed to be i.i.d. standard normal random variables, giving probit mediator and outcome models. 
This gives
\begin{align}
&\mathbb{E}\left( M_i|Z_i=z,\boldsymbol{X}_i=\boldsymbol{x}\right) = P\left( M_i=1|Z_i=z,\boldsymbol{X}_i=\boldsymbol{x}\right)   = \Phi\left( \beta_0 + \beta_1z + \boldsymbol{\beta}_2^\top\boldsymbol{x} + \boldsymbol{\beta}_3^\top z\boldsymbol{x}\right) , \label{eq:expMprobit} \\
&\mathbb{E}\left( Y_i|Z_i=z,M_i=m,\boldsymbol{X}_i=\boldsymbol{x}\right) =\Phi\Big( \theta_0+\theta_1z + \theta_2m + \theta_3zm + \boldsymbol{\theta}_4^\top\boldsymbol{x} + \boldsymbol{\theta}_5^\top z\boldsymbol{x} + \boldsymbol{\theta}_6^\top m\boldsymbol{x} + \boldsymbol{\theta}_7^\top zm\boldsymbol{x}\Big) , \label{eq:expYprobit}
\end{align} 
where $\Phi\left( \cdot\right) $ is the standard normal CDF. Substituting these into \eqref{eq:thm1nde} and \eqref{eq:thm1nie} yields expressions for the conditional natural direct and indirect effects
\begin{align}
&NDE(\boldsymbol{x}) =\left\{\Phi\left(\theta_0+\theta_1+\left( \boldsymbol{\theta}_4^\top+\boldsymbol{\theta}_5^\top\right) \boldsymbol{x}\right) - \Phi\left(\theta_0+\boldsymbol{\theta}_4^\top\boldsymbol{x}\right)\right\}\left(1-\Phi\left(\beta_0+\boldsymbol{\beta}_2^\top\boldsymbol{x}\right)\right) + \notag\\
&\quad\left\{\Phi\left(\theta_0+\theta_1+\theta_2+\theta_3+\left(\boldsymbol{\theta}_4^\top+\boldsymbol{\theta}_5^\top+\boldsymbol{\theta}_6^\top+\boldsymbol{\theta}_7^\top\right)\boldsymbol{x}\right) -\Phi\left(\theta_0+\theta_2+\left(\boldsymbol{\theta}_4^\top+\boldsymbol{\theta}_6^\top\right)\boldsymbol{x}\right)\right\}\Phi\left(\beta_0+\boldsymbol{\beta}_2^\top\boldsymbol{x}\right), \label{NDEhat} \\[7pt] 
&NIE(\boldsymbol{x}) =\left\{\Phi\left(\theta_0+\theta_1+\theta_2+\theta_3+\left(\boldsymbol{\theta}_4^\top+\boldsymbol{\theta}_5^\top+\boldsymbol{\theta}_6^\top+\boldsymbol{\theta}_7^\top\right)\boldsymbol{x}\right) - \Phi\left(\theta_0+\theta_1+\left( \boldsymbol{\theta}_4^\top+\boldsymbol{\theta}_5^\top\right) \boldsymbol{x}\right)\right\}\notag\\ 
&\quad\times\left\{\Phi\left(\beta_0+\beta_1+\left(\boldsymbol{\beta}_2^\top+\boldsymbol{\beta}_3^\top\right)\boldsymbol{x}\right)- \Phi\left(\beta_0+\boldsymbol{\beta}_2^\top\boldsymbol{x}\right)\right\}. \label{NIEhat}
\end{align}
Estimation can be performed by fitting \eqref{mstar}-\eqref{ystar} by maximum likelihood (ML). As \eqref{NDEhat}-\eqref{NIEhat} are functions of ML estimators, inference for $NDE(\boldsymbol{x})$ and $NIE(\boldsymbol{x})$ can be based on the delta method. The marginal (population averaged) effects $NDE$ and $NIE$ can be estimated by averaging the estimated conditional effects over the study population or sample, e.g. $\widehat{NIE}=\frac{1}{n}\sum_{i=1}^n\widehat{NIE}(\boldsymbol{x}_i)$, where $n$ is the size of the study population and  $\boldsymbol{x}_i $ is the covariate vector that has been observed for patient $i$.
Expressions for the alternative definitions of the natural direct and indirect effects can be similarly obtained (see Appendix A).

\section{Sensitivity analysis}
\label{sec:sensanalys}

\begin{figure}[ht]
	\centering	
	\includegraphics[]{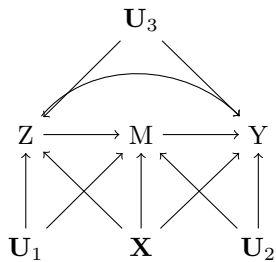}
	\caption{A directed acyclic graph with an exposure $Z$, a mediator $M$, an outcome $Y$, the set of observed confounders $\boldsymbol{X}$, and the unobserved confounders $\mathbf{U}_1$, $\mathbf{U}_2$, and $\mathbf{U}_3$.}
	\label{fig:DAG3}
\end{figure}

Figure \ref{fig:DAG3} illustrates three types of confounding relevant to the mediation setting, exposure-mediator confounders ($\mathbf{U}_1$), mediator-outcome confounders ($\mathbf{U}_2$), and exposure-outcome confounders ($\mathbf{U}_3$). The sensitivity analyses introduced in the literature consider the possible existence of $\mathbf{U}_2$. Here we consider all three potential sources of unobserved confounding $\mathbf{U}_1$, $\mathbf{U}_2$ and $\mathbf{U}_3$.

The techniques used here were first introduced by Genb\"{a}ck et al. \cite{genback2014} in the context of sensitivity analysis for linear regression parameters in the presence of non-ignorable missingness on a continuous outcome, and later generalized to binary outcomes \cite{genback2016}. 

\subsection{Uncertainty intervals for unobserved exposure-mediator confounding}
\label{sec:uizm}
Suppose that we can model the mediator $M_i$ as a function of the exposure and observed covariates as in \eqref{mstar} and that we can model the exposure assignment mechanism as a function of the observed covariates as $Z_i = I(Z_i^*>0)$, with
	\begin{equation}
	Z_i^*=\alpha_0 + \boldsymbol{\alpha}_1^\top\boldsymbol{X}_i + \varepsilon_i, \label{eq:zstar}
	\end{equation}
	where $\varepsilon_i$ are i.i.d. standard normal variables. If there is unobserved mediator-outcome confounding ($\mathbf{U}_1$ in Figure \ref{fig:DAG3}) this will induce a correlation between the error terms in the models for the exposure and the mediator, a fact that we will use in our sensitivity analysis.  
	
	Suppose that $\varepsilon_i$ and $\eta_{i}$ (the error term in \eqref{mstar}) are jointly normal with correlation $\rho_{\varepsilon\eta}$. If part \eqref{eg:seqign1} of Assumption \ref{ass:seqign} is fulfilled, i.e. the exposure-mediator relationship is unconfounded given the observed covariates $\boldsymbol{X}_i$, then $\rho_{\varepsilon\eta}=0$, otherwise $\rho_{\varepsilon\eta}\neq0$. If we take $\rho_{\varepsilon\eta}$ into account in the estimation of the regression parameters, which are then used in \eqref{NDEhat} and \eqref{NIEhat} to obtain the estimated $NDE(\boldsymbol{x})$ and $NIE(\boldsymbol{x})$, we get an idea of the effect of unobserved mediator-outcome confounding on the estimated conditional natural direct and indirect effects. 
	
	Let us denote the vectors of regression parameters in \eqref{mstar} and \eqref{eq:zstar} as  $\boldsymbol{\alpha}$, and $\boldsymbol{\beta}$. We can derive the log-likelihood of the regression parameters in \eqref{mstar} and \eqref{eq:zstar} and the correlation $\rho_{\varepsilon\eta}$, given the observed data, as  
	\begin{equation}
	\ell\left( \boldsymbol{\alpha},\boldsymbol{\beta},\rho_{\varepsilon\eta}\right)=\sum_{i}\left( 1-z_i\right)\ln\left\lbrace \Phi_2\left( w_{1i},-\boldsymbol{\alpha}^\top\boldsymbol{x}_i;-\rho_{1i}^*\right) \right\rbrace + \sum_{i}z_i\ln\left\lbrace \Phi_2\left( w_{1i},\boldsymbol{\alpha}^\top\boldsymbol{x}_i;\rho_{1i}^*\right) \right\rbrace, \label{eq:loglikm1}
	\end{equation}   
	where $ \Phi_2\left(\cdot ,\cdot;\cdot\right)$ is the standard bivariate normal cdf with three arguments, the first two are the means of the two random variables and the third is their correlation. We also have that $w_{1i}=\left(2m_i-1 \right)\left( \beta_0 + \beta_1z_i + \boldsymbol{\beta}_2^\top\boldsymbol{x}_i + \boldsymbol{\beta}_3^\top z_i\boldsymbol{x}_i\right)$, and $\rho_{1i}^*=\left(2m_i-1 \right)\rho_{\varepsilon\eta}$ \cite{genback2016,greene1998}. Using a modified maximum likelihood (ML) procedure where \eqref{eq:loglikm1} is maximized with regards to $\boldsymbol{\beta}$ and $\boldsymbol{\alpha}$ for a fixed $\rho_{\varepsilon\eta}=\tilde{\rho}_{\varepsilon\eta}$ we obtain $\hat{\boldsymbol{\beta}} \left( \tilde{\rho}_{\varepsilon\eta}\right) $, the estimated regression parameters in model \eqref{mstar} under correlation $\tilde{\rho}_{\varepsilon\eta}$. 

The estimate $\hat{\boldsymbol{\beta}} \left( \tilde{\rho}_{\varepsilon\eta}\right) $ can in turn be used (together with the $\hat{\boldsymbol{\theta}}$ obtained by fitting \eqref{ystar})
in \eqref{NDEhat} and \eqref{NIEhat} to obtain $\widehat{NDE}(\boldsymbol{x},\tilde{\rho}_{\varepsilon\eta})$ and $\widehat{NIE}(\boldsymbol{x},\tilde{\rho}_{\varepsilon\eta})$, estimates of the conditional natural direct and indirect effects under a given level of exposure-mediator confounding. 

The resulting $\widehat{NDE}(\boldsymbol{x},\tilde{\rho}_{\varepsilon\eta})$ and $\widehat{NIE}(\boldsymbol{x},\tilde{\rho}_{\varepsilon\eta})$ can be reported in different ways. One alternative is to plot these together with their confidence intervals over an interval of correlations (this is exemplified in Section \ref{sec:ressens}). Another alternative is to report the results through estimated \emph{identification sets}, consisting of the lower and upper bounds of the $\widehat{NDE}(\boldsymbol{x},\tilde{\rho}_{\varepsilon\eta})$ and $\widehat{NIE}(\boldsymbol{x},\tilde{\rho}_{\varepsilon\eta})$ over an interval of correlations. The estimated identification sets for the $NDE(\boldsymbol{x})$ and $NIE(\boldsymbol{x})$ over the interval $\tilde{\rho}_{\varepsilon\eta}\in\left[ a,b\right]$ are thus given by 
$$\left( \min_{\tilde{\rho}_{\varepsilon\eta}\in\left[ a,b\right]}\widehat{NDE}\left( \boldsymbol{x},\tilde{\rho}_{\varepsilon\eta}\right); \max_{\tilde{\rho}_{\varepsilon\eta}\in\left[ a,b\right]}\widehat{NDE}\left( \boldsymbol{x},\tilde{\rho}_{\varepsilon\eta}\right)\right),$$ and $$\left( \min_{\tilde{\rho}_{\varepsilon\eta}\in\left[ a,b\right]}\widehat{NIE}\left( \boldsymbol{x},\tilde{\rho}_{\varepsilon\eta}\right); \max_{\tilde{\rho}_{\varepsilon\eta}\in\left[ a,b\right]}\widehat{NIE}\left( \boldsymbol{x},\tilde{\rho}_{\varepsilon\eta}\right)\right).$$

To incorporate sampling variability \emph{uncertainty intervals} (UIs) \cite{vansteelandt2006} are constructed by taking the union of all confidence intervals obtained for the $NDE(\boldsymbol{x})$ and $NIE(\boldsymbol{x})$ with the correlation $\tilde{\rho}_{\varepsilon\eta}$ varying in the interval $\left[ a,b\right]$. The standard errors for $\widehat{NDE}\left( \boldsymbol{x},\tilde{\rho}_{\varepsilon\eta}\right)$ and $\widehat{NIE}\left( \boldsymbol{x},\tilde{\rho}_{\varepsilon\eta}\right)$ (see Appendix B) are used to construct $\left( 1-\alpha\right)\times100\% $ confidence intervals for $NDE(\boldsymbol{x})$ and $NIE(\boldsymbol{x})$. Let $LCI^{NDE}\left( \boldsymbol{x},\tilde{\rho}_{\varepsilon\eta}\right)$ ($LCI^{NIE}\left( \boldsymbol{x},\tilde{\rho}_{\varepsilon\eta}\right)$) and $UCI^{NDE}\left( \boldsymbol{x},\tilde{\rho}_{\varepsilon\eta}\right)$ ($UCI^{NIE}\left( \boldsymbol{x},\tilde{\rho}_{\varepsilon\eta}\right)$) denote the lower and upper bounds for the $\left( 1-\alpha\right)\times100\% $ CI of $NDE(\boldsymbol{x})$ ($NIE(\boldsymbol{x})$) for $\rho_{\varepsilon\eta}=\tilde{\rho}_{\varepsilon\eta}$. The lower and upper bounds of the (at least) $\left( 1-\alpha\right)\times100\% $ UI for exposure-mediator confounding are then given by 
$$ NDE(\boldsymbol{x})_{l,\rho_{\varepsilon\eta}}=\min_{\tilde{\rho}_{\varepsilon\eta}\in\left[ a,b\right]}LCI^{NDE}\left( \boldsymbol{x},\tilde{\rho}_{\varepsilon\eta}\right); NDE(\boldsymbol{x})_{u,\rho_{\varepsilon\eta}}=\max_{\tilde{\rho}_{\varepsilon\eta}\in\left[ a,b\right]}UCI^{NDE}\left( \boldsymbol{x},\tilde{\rho}_{\varepsilon\eta}\right),$$ and $$NIE(\boldsymbol{x})_{l,\rho_{\varepsilon\eta}}= \min_{\tilde{\rho}_{\varepsilon\eta}\in\left[ a,b\right]}LCI^{NIE}\left( \boldsymbol{x},\tilde{\rho}_{\varepsilon\eta}\right); NIE(\boldsymbol{x})_{u,\rho_{\varepsilon\eta}}=\max_{\tilde{\rho}_{\varepsilon\eta}\in\left[ a,b\right]}UCI^{NIE}\left( \boldsymbol{x},\tilde{\rho}_{\varepsilon\eta}\right).$$

A plausible interval of correlations could be narrowed down using subject-matter knowledge. For example, is it more likely that the correlation induced by an unobserved confounder is positive or negative? To understand the connection between the size of the correlation and strength of confounding one could, e.g., investigate the observed covariates and the effect of leaving out the strongest confounder.

It is often of interest to ascertain the degree of unobserved confounding that would render an effect non-significant. Therefore, in addition to reporting the UIs themselves (or as an alternative) one could report ranges of correlations where the $\left( 1-\alpha\right)\times100\% $ UI includes 0 (i.e. where the effect is not significant at the $\alpha$ significance level). This approach is also exemplified in Section \ref{sec:ressens}.

A sensitivity analysis of the marginal effects can be performed by averaging the $\widehat{NDE}(\boldsymbol{x},\tilde{\rho}_{\varepsilon\eta})$ ($\widehat{NIE}(\boldsymbol{x},\tilde{\rho}_{\varepsilon\eta})$) over the study population to obtain $\widehat{NDE}(\tilde{\rho}_{\varepsilon\eta})$ ($\widehat{NIE}(\tilde{\rho}_{\varepsilon\eta})$). The corresponding standard errors can be obtained using the delta method (see Appendix B) and estimated identification sets and UIs for the marginal effects under exposure-mediator confounding constructed as outlined above.

\subsection{Uncertainty intervals for unobserved mediator-outcome confounding}
\label{sec:uimy}
We now turn our attention to part \eqref{eg:seqign2} of Assumption \ref{ass:seqign} which states that, given the observed exposure $Z_i$ and the observed covariates $\boldsymbol{X}_i$ the mediator-outcome relation is unconfounded. Suppose that the observed mediator can be modeled as in \eqref{mstar} and that the observed outcome can be modeled as in \eqref{ystar}.
We assume that $\eta$ and $\xi$ (the error terms in \eqref{mstar} and \eqref{ystar}) are bivariate standard normal distributed with correlation $\rho_{\eta\xi}$. If part \eqref{eg:seqign2} of Assumption \ref{ass:seqign} is fulfilled then $\rho_{\eta\xi}=0$, otherwise $\rho_{\eta\xi}\neq0$.
	
We denote the vector of regression parameters in \eqref{ystar} as  $\boldsymbol{\theta}$. The log-likelihood of the regression parameters in \eqref{mstar} and \eqref{ystar} and the correlation $\rho_{\eta\xi}$, given the observed data, is given by
\begin{equation}
\hspace{-0.05cm}\ell\left( \boldsymbol{\beta},\boldsymbol{\theta},\rho_{\eta\xi}\right)=\sum_{i}\left( 1-m_i\right)\ln\left\lbrace \Phi_2\left( w_{2i},-\boldsymbol{\beta}^\top\boldsymbol{c}_{i};-\rho_{2i}^*\right) \right\rbrace + \sum_{i}m_i\ln\left\lbrace \Phi_2\left( w_{2i},\boldsymbol{\beta}^\top\boldsymbol{c}_{i};\rho_{2i}^*\right) \right\rbrace. \label{eq:logliky1}
\end{equation}
Here $w_{2i}=\left( 2y_i-1\right)\left( \theta_0+\theta_1z_i + \theta_2m_i + \theta_3z_im_i + \boldsymbol{\theta}_4^\top\boldsymbol{x}_i + \boldsymbol{\theta}_5^\top z_i\boldsymbol{x}_i + \boldsymbol{\theta}_6^\top m_i\boldsymbol{x}_i + \boldsymbol{\theta}_7^\top z_im_i\boldsymbol{x}_i\right) $, $\boldsymbol{c}_{i}=(z_i,\boldsymbol{x}_i^\top,z_i\boldsymbol{x}_i^\top)^\top$ and $\rho_{2i}^*=\left(2y_i-1 \right)\rho_{\eta\xi}$. By maximizing \eqref{eq:logliky1} with regards to $\boldsymbol{\theta}$ and $\boldsymbol{\beta}$ for a fixed $\rho_{\eta\xi}=\tilde{\rho}_{\eta\xi}$ we obtain $\hat{\boldsymbol{\theta}} \left( \tilde{\rho}_{\eta\xi}\right) $ and $\hat{\boldsymbol{\beta}}\left( \tilde{\rho}_{\eta\xi}\right)$, the estimated regression parameters in models \eqref{ystar} and \eqref{mstar} under correlation $\tilde{\rho}_{\eta\xi}$.  
	
The $\hat{\boldsymbol{\theta}}\left( \tilde{\rho}_{\eta\xi}\right) $ and  $\hat{\boldsymbol{\beta}}\left( \tilde{\rho}_{\eta\xi}\right)$, allow us to calculate $\widehat{NDE}(\boldsymbol{x},\tilde{\rho}_{\eta\xi})$  and $\widehat{NIE}(\boldsymbol{x},\tilde{\rho}_{\eta\xi})$.
Estimated identification sets for the $NDE(\boldsymbol{x})$ and $NIE(\boldsymbol{x})$ for $\tilde{\rho}_{\eta\xi}\in\left[ a',b'\right] $ are then given by $$\left( \min_{\tilde{\rho}_{\eta\xi}\in\left[ a',b'\right]}\widehat{NDE}\left( \boldsymbol{x},\tilde{\rho}_{\eta\xi}\right); \max_{\tilde{\rho}_{\eta\xi}\in\left[ a',b'\right]}\widehat{NDE}\left( \boldsymbol{x},\tilde{\rho}_{\eta\xi}\right)\right),$$ and $$\left( \min_{\tilde{\rho}_{\eta\xi}\in\left[ a',b'\right]}\widehat{NIE}\left( \boldsymbol{x},\tilde{\rho}_{\eta\xi}\right); \max_{\tilde{\rho}_{\eta\xi}\in\left[ a',b'\right]}\widehat{NIE}\left( \boldsymbol{x},\tilde{\rho}_{\eta\xi}\right)\right).$$
Let $LCI^{NDE}\left( \boldsymbol{x},\tilde{\rho}_{\eta\xi}\right)$ ($LCI^{NIE}\left( \boldsymbol{x},\tilde{\rho}_{\eta\xi}\right)$) and $UCI^{NDE}\left( \boldsymbol{x},\tilde{\rho}_{\eta\xi}\right)$ ($UCI^{NIE}\left( \boldsymbol{x},\tilde{\rho}_{\eta\xi}\right)$) denote the lower and upper bounds for the $\left( 1-\alpha\right)\times100\% $ CI of $NDE(\boldsymbol{x})$ ($NIE(\boldsymbol{x})$) for $\rho_{\eta\xi}=\tilde{\rho}_{\eta\xi}$. At least $\left( 1-\alpha\right)\times100\% $ UIs under mediator-outcome confounding are then given by the lower and upper bounds
$$ NDE(\boldsymbol{x})_{l,\rho_{\eta\xi}}=\min_{\tilde{\rho}_{\eta\xi}\in\left[ a',b'\right]}LCI^{NDE}\left( \boldsymbol{x},\tilde{\rho}_{\eta\xi}\right); NDE(\boldsymbol{x})_{u,\rho_{\eta\xi}}=\max_{\tilde{\rho}_{\eta\xi}\in\left[ a',b'\right]}UCI^{NDE}\left( \boldsymbol{x},\tilde{\rho}_{\eta\xi}\right),$$ and $$NIE(\boldsymbol{x})_{l,\rho_{\eta\xi}}= \min_{\tilde{\rho}_{\eta\xi}\in\left[ a',b'\right]}LCI^{NIE}\left( \boldsymbol{x},\tilde{\rho}_{\eta\xi}\right); NIE(\boldsymbol{x})_{u,\rho_{\eta\xi}}=\max_{\tilde{\rho}_{\eta\xi}\in\left[ a',b'\right]}UCI^{NIE}\left( \boldsymbol{x},\tilde{\rho}_{\eta\xi}\right).$$

The marginal effects $\widehat{NDE}(\tilde{\rho}_{\eta\xi})$ and $\widehat{NIE}(\tilde{\rho}_{\eta\xi})$ are obtained by averaging the $\widehat{NDE}(\boldsymbol{x},\tilde{\rho}_{\eta\xi})$ and $\widehat{NIE}(\boldsymbol{x},\tilde{\rho}_{\eta\xi})$ over the study population and standard errors obtained through the delta method. Estimated identification sets and UIs for the $NDE$ and $NIE$ under mediator-outcome confounding can then be constructed as outlined above.

\subsection{Uncertainty intervals for unobserved exposure-outcome confounding}
\label{sec:uizy}
	
Finally, we address the issue of unobserved exposure-outcome confounding (i.e. $\mathbf{U}_3$ in Figure \ref{fig:DAG3}). 
Suppose again that we can model the exposure as \eqref{eq:zstar} and the outcome as \eqref{ystar}. We assume that  $\varepsilon_i$ and $\xi_i$ (the error terms in \eqref{eq:zstar} and \eqref{ystar}) are bivariate standard normal distributed with correlation $\rho_{\varepsilon\xi}$. If there is no unobserved exposure-outcome confounding then $\rho_{\varepsilon\xi}=0$, otherwise $\rho_{\varepsilon\xi}\neq0$.
	
The log-likelihood of the regression parameters in \eqref{eq:zstar} and \eqref{ystar} and the correlation $\rho_{\varepsilon\xi}$, given the observed data, is given by
\begin{equation}
\ell\left( \boldsymbol{\alpha},\boldsymbol{\theta},\rho_{\varepsilon\xi}\right)=\sum_{i}\left( 1-z_i\right)\ln\left\lbrace \Phi_2\left( w_{2i},-\boldsymbol{\alpha}^\top\boldsymbol{x}_i;-\rho_{3i}^*\right) \right\rbrace + \sum_{i}z_i\ln\left\lbrace \Phi_2\left( w_{2i},\boldsymbol{\alpha}^\top\boldsymbol{x}_i;\rho_{3i}^*\right) \right\rbrace. \label{eq:loglikyz1}
\end{equation}  
Here $\rho_{3i}^*=\left(2y_i-1 \right)\rho_{\varepsilon\xi}$ and $w_{2i}$ as before. By maximizing (\ref{eq:loglikyz1}) with regards to $\boldsymbol{\theta}$ and $\boldsymbol{\alpha}$ for a fixed $\rho_{\varepsilon\xi}=\tilde{\rho}_{\varepsilon\xi}$ we obtain $\hat{\boldsymbol{\theta}} \left( \tilde{\rho}_{\varepsilon\xi}\right) $, the estimated regression parameters in model (\ref{ystar}) under correlation $\tilde{\rho}_{\varepsilon\xi}$.  
	
Using $\hat{\boldsymbol{\theta}}\left( \tilde{\rho}_{\varepsilon\xi}\right) $ and the $\hat{\boldsymbol{\beta}}$ obtained from fitting \eqref{mstar} in \eqref{NDEhat}-\eqref{NIEhat} gives us  $\widehat{NDE}(\boldsymbol{x},\tilde{\rho}_{\varepsilon\xi})$  and $\widehat{NIE}(\boldsymbol{x},\tilde{\rho}_{\varepsilon\xi})$.
Estimated identification sets for the $NDE(\boldsymbol{x})$ and $NIE(\boldsymbol{x})$ for $\tilde{\rho}_{\varepsilon\xi}\in\left[ a^*,b^*\right] $ are then given by $$\left( \min_{\tilde{\rho}_{\varepsilon\xi}\in\left[ a^*,b^*\right]}\widehat{NDE}\left( \boldsymbol{x},\tilde{\rho}_{\varepsilon\xi}\right); \max_{\tilde{\rho}_{\varepsilon\xi}\in\left[ a^*,b^*\right]}\widehat{NDE}\left( \boldsymbol{x},\tilde{\rho}_{\varepsilon\xi}\right)\right),$$ and $$\left( \min_{\tilde{\rho}_{\varepsilon\xi}\in\left[ a^*,b^*\right]}\widehat{NIE}\left( \boldsymbol{x},\tilde{\rho}_{\varepsilon\xi}\right); \max_{\tilde{\rho}_{\varepsilon\xi}\in\left[ a^*,b^*\right]}\widehat{NIE}\left( \boldsymbol{x},\tilde{\rho}_{\varepsilon\xi}\right)\right).$$
Thus, at least $\left( 1-\alpha\right)\times100\% $ UIs under mediator-outcome confounding are given by the lower and upper bounds
$$ NDE(\boldsymbol{x})_{l,\rho_{\varepsilon\xi}}=\min_{\tilde{\rho}_{\varepsilon\xi}\in\left[ a^*,b^*\right]}LCI^{NDE}\left( \boldsymbol{x},\tilde{\rho}_{\varepsilon\xi}\right); NDE(\boldsymbol{x})_{u,\rho_{\eta\xi}}=\max_{\tilde{\rho}_{\varepsilon\xi}\in\left[ a^*,b^*\right]}UCI^{NDE}\left( \boldsymbol{x},\tilde{\rho}_{\varepsilon\xi}\right),$$ and $$NIE(\boldsymbol{x})_{l,\rho_{\varepsilon\xi}}= \min_{\tilde{\rho}_{\varepsilon\xi}\in\left[ a^*,b^*\right]}LCI^{NIE}\left( \boldsymbol{x},\tilde{\rho}_{\varepsilon\xi}\right); NIE(\boldsymbol{x})_{u,\rho_{\varepsilon\xi}}=\max_{\tilde{\rho}_{\varepsilon\xi}\in\left[ a^*,b^*\right]}UCI^{NIE}\left( \boldsymbol{x},\tilde{\rho}_{\varepsilon\xi}\right).$$
	
Again, estimated identification sets and UIs for marginal effects under exposure-outcome confounding can be obtained by averaging the $\widehat{NDE}(\boldsymbol{x},\tilde{\rho}_{\varepsilon\xi})$ and $\widehat{NIE}(\boldsymbol{x},\tilde{\rho}_{\varepsilon\xi})$ over the study population and using the delta method for the corresponding standard errors. 

Finally, note that the suggested methods evaluate sensitivity to each type of unobserved confounding separately, assuming that the other two types are not present.

\section{Case study}
\label{sec:casestudy}
Previous studies have shown evidence that living alone is detrimental to the outcome after stroke, especially in male patients \cite{lindmark2014,eriksson2009}.  One possible explanation is that patients living alone are less likely to recognize stroke symptoms and therefore arrive to the hospital later and with a more severe stroke than patients who cohabit. With a focus on male patients, we have used data from Riksstroke, the Swedish Stroke Register, to investigate to what extent the effect of living alone on the outcome after stroke is mediated by stroke severity (Section \ref{sec:reseffects}). We used the methods proposed in Section  \ref{sec:sensanalys} to assess the sensitivity of our findings to unobserved confounding (Section \ref{sec:ressens}).

\subsection{Data}
Riksstroke was established in 1994 with the main purpose of monitoring and supporting quality improvement of the Swedish stroke care. It covers all Swedish hospitals that admit acute stroke patients and patient-level information is collected during the acute phase and at follow-up 3 months and 1 year after stroke \cite{asplund2011}.

The data used in this example consist of 7 639 men with intracerebral hemorrhage or ischemic stroke (called simply ``stroke" in the sequel) who were registered in Riksstroke between January 1 and October 1, 2012. Patients included were registered as living at home and being independent in activities of daily living (ADL) at the time of stroke. 

The binary outcome variable was death or dependency in ADL at 3 months, defined as the patient being registered as dependent in ADL at the 3 month follow-up or the patient dying within 90 days after their stroke. Dependence in ADL was defined as the patient being unable to manage dressing, using the toilet, or walking indoors unassisted. Dates of death were retrieved from the Swedish Cause of Death Register managed by the National Board of Health and Welfare. 

The exposure variable was whether or not the patient was living alone at the time of stroke and the mediator variable was whether the patient had lowered consciousness at admission to the hospital or was fully conscious. The level of consciousness at admission was used as a proxy for stroke severity and corresponds to two levels based on the Reaction Level Scale (RLS) \cite{starmark1988} where fully conscious corresponds to RLS 1, and lowered consciousness corresponds to RLS 2–-8. 

Adjustment for confounding was made using the available pre-exposure covariates: highest attained education level and patient age at the time of stroke. Highest attained education level was obtained from the LISA database (Longitudinal integration database for health insurance and labor market studies) managed by Statistics Sweden, and was categorized into two groups; patients with and without university education. Age was modeled using both a continuous and categorical variable to take into account effect differences in different age groups. These age groups were allowed to differ between the exposure, mediator, and outcome models, depending on the best fit for each model. 

We performed a complete case analysis, meaning that cases with missing data on either the outcome or any of the covariates were deleted. Thus, 1 016 cases were deleted prior to analysis due to missing outcomes.  These missing outcomes consisted of patients who survived the three month mark but who did not provide follow-up information on ADL dependency. The primary reason for missing follow-up data is likely to be that the hospital where the patient was treated has not sent out the follow-up questionnaire. This is unlikely to be correlated to the ADL-dependency of the patient, and thus it is plausible to regard the outcomes as being missing at random. The other variables had a smaller number of missing values (22-209 cases). The final data used for the analyses consisted of a total of 6 432 patients.  

Analyses were performed using the R software environment, version 3.1.0 \cite{r2015}.

\subsection{Results}
The patients included in this study were on average 72.2 years old (range 18-98 years), and 1 258 (19.6\%) had a university education. A total number of 2 058 (32.0\%) patients were living alone before their stroke and a larger proportion of these had a lowered consciousness level on arrival to hospital compared to cohabitant patients (14.4\% vs. 11.2\%). Patients living alone were also more often dead or dependent in ADL 3 months after stroke (28.8\% vs. 23.9\%).

\subsubsection{Mediation effects}
\label{sec:reseffects}
Tables \ref{table:expmod}-\ref{table:outcmod} show the estimated probit models that are the basis for the analyses. The estimate of the exposure model \eqref{eq:zstar} showed a negative association between university education and the probability of living alone at the time of stroke (Table \ref{table:expmod}). Although the effect of age as a continuous variable was negative, older age groups had a significantly higher probability of living alone compared to patients under the age of 80.

\begin{table}[ht]
	\caption{Estimated probit model for the exposure living alone. Estimated regression parameters and standard errors.}
	\begin{center}
		\begin{tabular}{l D{.}{.}{2.6} D{.}{.}{2.6}}
			\toprule
			& \multicolumn{1}{c}{Estimated parameter} & \multicolumn{1}{c}{Standard error}\\
			\midrule
			Intercept            & -0.154     & 0.138      \\
			\noalign{\vskip 1mm}
			Education: &  &\\
			\quad No university & \multicolumn{1}{c}{Ref.}  & \multicolumn{1}{c}{Ref.} \\
			\quad University & -0.242^{***} & 0.043      \\
			\noalign{\vskip 1mm}
			Age (continuous)     & -0.005^{**}  & 0.002      \\
			\noalign{\vskip 1mm}
			Age (categorical): &  &\\
			\quad 18-79      &   \multicolumn{1}{c}{Ref.} 	&  \multicolumn{1}{c}{Ref.}\\
			\quad 80-84            & 0.256^{***} 		& 0.056      \\
			\quad 85-89            & 0.378^{***}  	& 0.068      \\
			\quad 90-              & 0.775^{***}  	& 0.090      \\
			\bottomrule
			\multicolumn{3}{l}{\scriptsize{$^{***}p<0.001$, $^{**}p<0.01$, $^*p<0.05$}}\\
		\end{tabular}
		\label{table:expmod}
	\end{center}
\end{table} 

In the estimated mediator model \eqref{mstar} there was a significantly higher probability of having lowered consciousness upon arrival to hospital for patients living alone compared to cohabitant patients, and older age groups compared to patients under the age of 75. There was also a significant interaction between living alone and having a university education (Table \ref{table:medmod}).

\begin{table}[ht]
	\caption{Estimated probit model for the mediator lowered consciousness. Estimated regression parameters and standard errors.}
	\begin{center}
		\begin{tabular}{l D{.}{.}{2.6} D{.}{.}{2.6} }
			\toprule
			& \multicolumn{1}{c}{Estimated parameter} & \multicolumn{1}{c}{Standard error} \\
			\midrule
			Intercept      & -1.169^{***} & 0.200\\
			\noalign{\vskip 1mm} 
			Cohabitation status:                     &       &       \\
			\quad Cohabitant                  &   \multicolumn{1}{c}{Ref.}      & \multicolumn{1}{c}{Ref.}   \\
			\quad Living alone                & 0.159^{***}    & 0.047   \\
			\noalign{\vskip 1mm} 
			Education:                     &       &       \\
			\quad No university         &   \multicolumn{1}{c}{Ref.}      & \multicolumn{1}{c}{Ref.}   \\
			\quad University            & 0.065     & 0.061   \\
			\noalign{\vskip 1mm} 
			Age (continuous)      		& -0.003     & 0.003   \\
			\noalign{\vskip 1mm} 
			Age (categorical):                     &       &   \\
			\quad 18-74                                &   \multicolumn{1}{c}{Ref.}      & \multicolumn{1}{c}{Ref.}    \\
			\quad 75-84                 & 0.205^{**} & 0.068  \\
			\quad 85-                   & 0.559^{***} & 0.095 \\
			\noalign{\vskip 1mm} 
			Living alone$\times$University &       -0.239^{*} & 0.118   \\
			\bottomrule
			\multicolumn{3}{l}{\scriptsize{$^{***}p<0.001$, $^{**}p<0.01$, $^*p<0.05$}}\\
		\end{tabular}
		\label{table:medmod}
	\end{center}
\end{table}

In the estimated outcome model \eqref{ystar} there was a significant positive effect of the mediator, level of consciousness (Table \ref{table:outcmod}). There was also a significant positive effect of the treatment, living alone, as well as a significant interaction between living alone and the age group 80-89. Age as a continuous variable and the older age groups, ages 80-89 and 90 and above, compared to the youngest age group, ages 79 and below, were positively associated with the probability of death or being dependent in ADL at three months after stroke. Finally, having a university education was significantly associated with a lowered probability of death or ADL dependency at three months.

\begin{table}[ht]
	\caption{Estimated probit model for the outcome dead or dependent in ADL at 3 months. Estimated regression parameters and standard errors.}
	\begin{center}
		\begin{tabular}{l D{.}{.}{2.6} D{.}{.}{2.6} }
			\toprule
			& \multicolumn{1}{c}{Estimated parameter} & \multicolumn{1}{c}{Standard error}  \\
			\midrule
			Intercept     & -2.764^{***} & 0.187\\
			\noalign{\vskip 1mm} 
			Cohabitation status:                     &       &       \\
			\quad Cohabitant                  &   \multicolumn{1}{c}{Ref.}      & \multicolumn{1}{c}{Ref.}    \\
			\quad Living alone     & 0.138^{**}   &  0.051   \\
			\noalign{\vskip 1mm} 
			Level of consciousness:                     &       &           \\
			\quad Fully cons.                  &   \multicolumn{1}{c}{Ref.}      & \multicolumn{1}{c}{Ref.}   \\
			\quad Lowered cons.        & 1.502^{***}     &  0.054 \\
			\noalign{\vskip 1mm} 
			Education:                     &       &        \\
			\quad No university                  &   \multicolumn{1}{c}{Ref.}      & \multicolumn{1}{c}{Ref.}   \\
			\quad University        & -0.115^{*}  & 0.049 \\
			\noalign{\vskip 1mm} 
			Age (continuous)     & 0.024^{***} & 0.003 \\
			\noalign{\vskip 1mm} 
			Age (categorical):                     &       &       \\
			\quad 18-79               & \multicolumn{1}{c}{Ref.}       & \multicolumn{1}{c}{Ref.}  \\
			\quad 80-89         & 0.397^{***}  & 0.065\\
			\quad 90-          & 0.619^{***}  & 0.132 \\
			\noalign{\vskip 1mm}
			Living alone$\times$80-89       & -0.300^{***}  & 0.086    \\
			Living alone$\times$90-         & -0.249  & 0.161      \\
			\noalign{\vskip 1mm}
			\bottomrule
			\multicolumn{3}{l}{\scriptsize{$^{***}p<0.001$, $^{**}p<0.01$, $^*p<0.05$}}\\
		\end{tabular} 
		\label{table:outcmod}
	\end{center}
\end{table}
		
We estimated the marginal $NDE$ and $NIE$ as well as the $NDE(\boldsymbol{x})$ and $NIE(\boldsymbol{x})$ for different covariate patterns, corresponding to a patient of average age (72.2 years old), average age minus one standard deviation (60.4), average age plus one standard deviation (84.1), and conditioning on level of education (Table \ref{table:effects}).

		\begin{table}[ht]
	\caption{Estimated marginal and conditional natural direct and indirect effects and total effects. 95\% CIs in parentheses.}
	\begin{center}
		\begin{tabular}{l D{.}{.}{-2.9} D{.}{.}{-2.9} D{.}{.}{-2.9}}
			\toprule
			\noalign{\vskip 1mm} 
			& \multicolumn{1}{c}{Natural direct effect} & \multicolumn{1}{c}{Natural indirect effect} & \multicolumn{1}{c}{Total effect} \\
			\midrule
			\noalign{\vskip 1mm}
			Marginal    & 0.006       & 0.012^{**}      &0.018 \\
			& \multicolumn{1}{c}{$(-0.014, 0.027)$}      & \multicolumn{1}{c}{$(0.003, 0.021)$}   &  \multicolumn{1}{c}{$(-0.003, 0.040)$} \\
			\noalign{\vskip 1mm}
			Conditional   &        &    &    \\
			\noalign{\vskip 1mm}
			\hspace{2mm}University education    &        &     &   \\
			\noalign{\vskip 1mm}
			\hspace{4mm}60.4 years     & 0.024^{**}       & -0.006     &  0.018 \\
			& \multicolumn{1}{c}{$(0.006,0.043 )$}      & \multicolumn{1}{c}{$(-0.023, 0.010)$}   &  \multicolumn{1}{c}{$(-0.005, 0.041)$} \\
			\noalign{\vskip 1mm}
			\hspace{4mm}72.2 years     & 0.032^{**}       & -0.007     & 0.025  \\
			& \multicolumn{1}{c}{$(0.008, 0.056)$}      & \multicolumn{1}{c}{$(-0.024,0.011)$}   & \multicolumn{1}{c}{$(-0.003,0.053)$}  \\
			\noalign{\vskip 1mm}
			\hspace{4mm}84.1 years     & -0.053^{*}       & -0.009   &  -0.062^{*}   \\
			& \multicolumn{1}{c}{$(-0.097, -0.009)$}      & \multicolumn{1}{c}{$(-0.031, 0.014)$}  &  \multicolumn{1}{c}{$(-0.112, -0.011)$}  \\
			\noalign{\vskip 1mm}
			\hspace{2mm}No university education    &        &   &     \\
			\hspace{4mm}60.4 years            & 0.027^{**} & 0.014^{**} & 0.041^{***}\\
			& \multicolumn{1}{c}{$(0.007, 0.047)$}      & \multicolumn{1}{c}{$(0.006, 0.023)$}   & \multicolumn{1}{c}{$(0.020, 0.063)$} \\
			\noalign{\vskip 1mm}
			\hspace{4mm}72.2 years            & 0.035^{**} & 0.015^{**} & 0.050^{***}\\
			& \multicolumn{1}{c}{$(0.009, 0.061)$}      & \multicolumn{1}{c}{$(0.006, 0.024)$}    & \multicolumn{1}{c}{$(0.023, 0.076)$} \\
			\noalign{\vskip 1mm}
			\hspace{4mm}84.1 years            & -0.056^{*} & 0.018^{**} & -0.038\\
			& \multicolumn{1}{c}{$(-0.102, -0.009)$}      & \multicolumn{1}{c}{$(0.007, 0.029)$}    & \multicolumn{1}{c}{$(-0.086, 0.011)$} \\
			\noalign{\vskip 1mm}
			\bottomrule
			\multicolumn{4}{l}{\scriptsize{$^{***}p<0.001$, $^{**}p<0.01$, $^*p<0.05$}}
		\end{tabular} 
		\label{table:effects}
	\end{center}
\end{table}
			
The marginal natural direct and indirect effects were both positive, with only the natural indirect effect significant. The marginal total effect of living alone on death or dependency in ADL was not significant at the 5\% level. All the conditional total effects were significant at the 5\% level except for those cases (University education and ages 60.4 or 72.2, No university education and age 84.1) where the natural direct and indirect effects had opposite signs. As indicated by the interaction between the exposure living alone and education level in the mediator model (Table \ref{table:medmod}) the conditional natural indirect effect differs substantially between patients with and without university education (second column of Table \ref{table:effects}), where the effects were negative and not significant for the former but positive and significant at the 5\% level for the latter. That is, for patients with university education and the three investigated ages there was no evidence of an indirect effect of living alone on death or ADL-dependency at 3 months working through level of consciousness. For patients without university education the positive significant indirect effects indicated that living alone increases the probability of having lowered consciousness upon arrival to hospital which in turn increases the probability of death or dependency in ADL at 3 months.
				
As suggested by the interaction between living alone and age group in the outcome model (Table \ref{table:outcmod}) the conditional natural direct effect (first column of Table \ref{table:effects}), i.e. the effect of living alone on death or ADL-dependency not working through differences in level of consciousness upon arrival, differs quite a bit between an 84.1 year old patient (falling in the 80-89 year category) and a 60.4 or 72.2 year old patient (both falling in the 18-79 year category). All conditional natural direct effects were significant at the 5\% level, but positive for 60.4 and 72.2 year old patients, meaning that living alone increases the probability of death or dependency in ADL, and negative for an 84.1 year old patient, indicating a decreased probability of death or dependency in ADL for patients living alone not through differences in level of consciousness. 
			
\subsubsection{Sensitivity analysis}
\label{sec:ressens}
We continue by investigating how sensitive the significant effects in Table \ref{table:effects} are to unobserved confounding. In our example possible unobserved exposure-mediator confounders ($\mathbf{U}_1$) include geographical factors such as distance from the patient's home to the hospital and possible unobserved mediator-outcome confounders ($\mathbf{U}_2$) include genetic factors. It is also possible that pre-exposure socioeconomic factors not captured by education level could confound the exposure-outcome relation ($\mathbf{U}_3$).
			
\begin{figure}[ht]
	\centering	
	\includegraphics[scale=0.75]{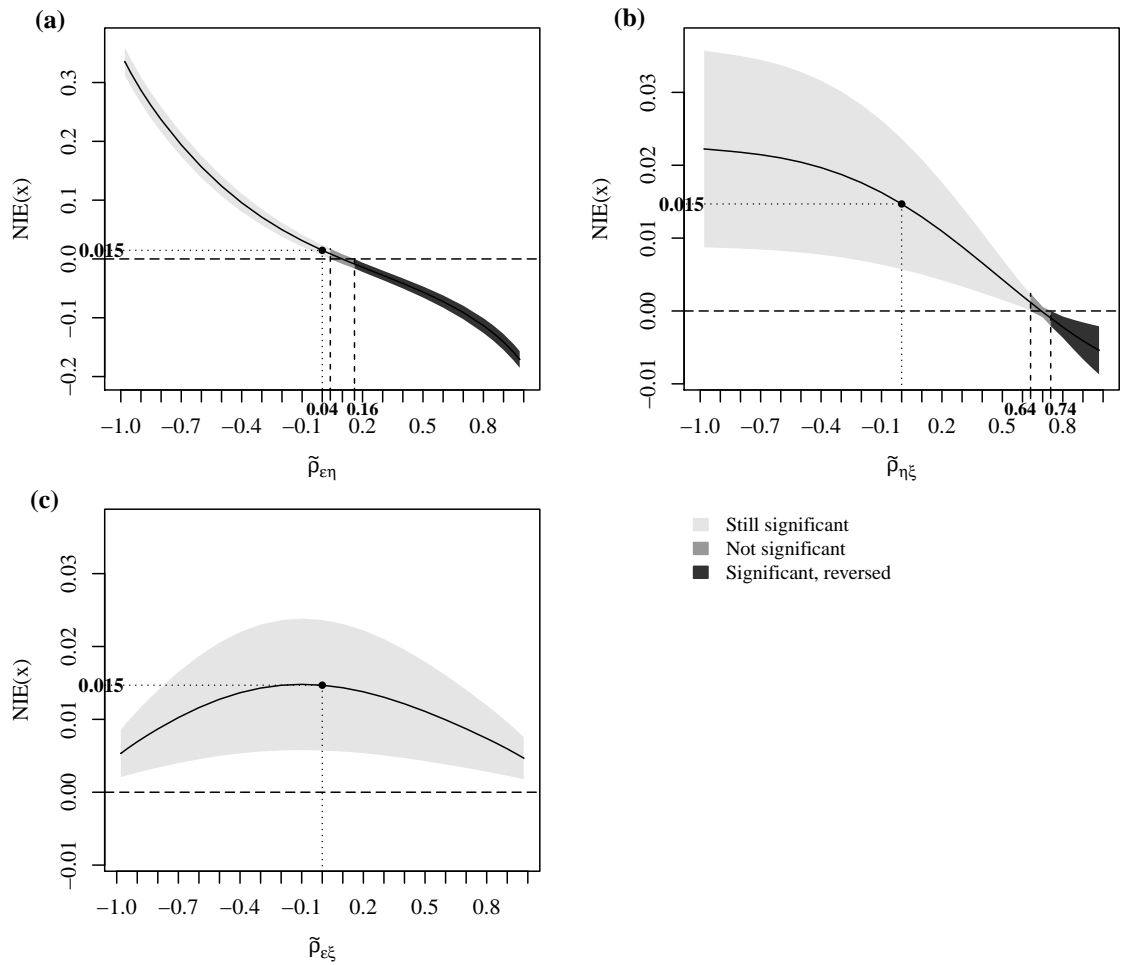}
	\caption{Estimated $ NIE(\boldsymbol{x}) $ for a patient of average age (72.2 years) without university education with corresponding 95\% CIs (shaded area) in the presence of \textbf{(a)} exposure-mediator, \textbf{(b)} mediator-outcome, and \textbf{(c)} exposure-outcome confounding. The light gray areas correspond to 95\% CIs that lie entirely above 0, the medium gray areas to 95\% CIs that include 0 and the dark gray areas to 95\% CIs where the effect is reversed. Note that the scale of the y-axis differs between panel \textbf{(a)} and panels \textbf{(b)}-\textbf{(c)}.}
	\label{fig:sensanMean}
\end{figure}
			
\begin{figure}[ht]
	\centering	
	\includegraphics[scale=0.75]{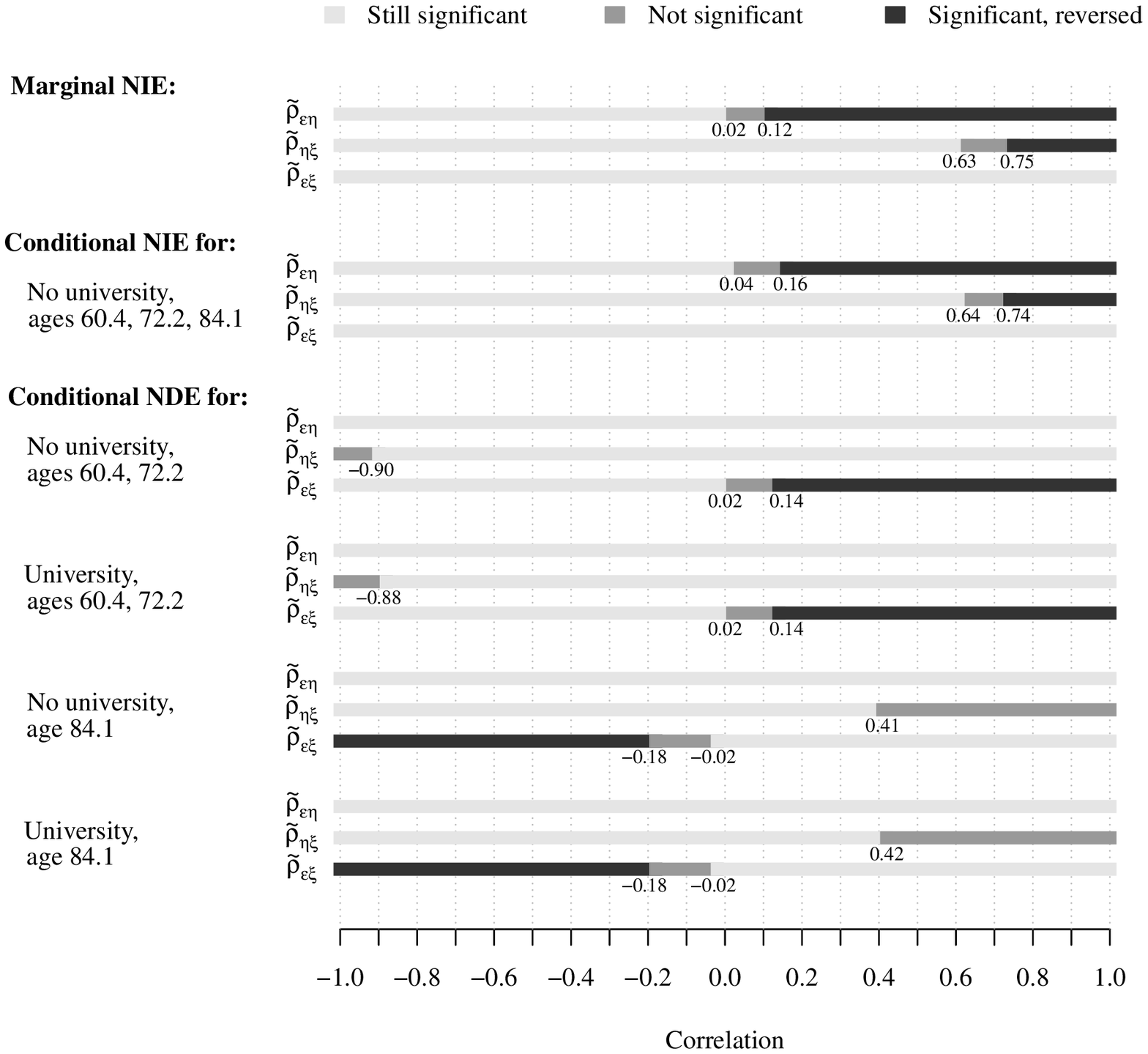}
	\caption{Results of the sensitivity analysis of the significant effects from Table \ref{table:effects}. Ranges of correlations ($\rho_{\varepsilon\eta}$: exposure-mediator confounding, $\rho_{\eta\xi}$: mediator-outcome confounding, and $\rho_{\varepsilon\xi}$: exposure-outcome confounding) that would render the effect significant and in the same direction as in Table \ref{table:effects} (light gray), not significant (medium gray) or reversed (dark gray) at the 5\% level.}
	\label{fig:summary}
\end{figure}
			
Figure \ref{fig:sensanMean} shows the estimated $NIE(\boldsymbol{x})$ with corresponding 95\% CIs for a 72.2 year-old patient with no university education under varying levels of exposure-mediator confounding (Figure \ref{fig:sensanMean}a), mediator-outcome confounding (Figure \ref{fig:sensanMean}b), and exposure-outcome confounding (Figure \ref{fig:sensanMean}c). The medium gray shaded areas in Figure \ref{fig:sensanMean} correspond to the 95\% CIs that include 0, i.e. where there is a non-significant effect, while the dark gray shaded areas correspond to 95\% CIs that lie entirely below 0, indicating a reversal of the effect. The light gray shaded areas correspond to 95\% CIs that lie entirely above 0, i.e. areas where the effect is still positive and significant at the 5\% level. The 95\% UI is given by the union of all 95\% CIs over a given interval of correlations.  From Figure \ref{fig:sensanMean}a we see that the $\widehat{NIE}(\boldsymbol{x})$ is largest for a correlation of $-0.98$ and then decreases as the correlation goes towards 1. We also see that $0 \notin 95\%$ UI for $\rho_{\varepsilon\eta}\in\left(-1,0.04 \right)$, and that the upper bound of the 95\% UI, $NIE(\boldsymbol{x})_{u,\rho_{\varepsilon\eta}}<0$, i.e. the effect is reversed, if $\rho_{\varepsilon\eta}\in\left(0.16,1 \right)$.  This indicates that the $\widehat{NIE}(\boldsymbol{x})$ is insensitive to exposure-mediator confounding that induces a negative correlation between the error terms in the exposure and mediator models, but quite sensitive to even a moderate positive correlation. A situation that would induce positive correlation between the error terms in the exposure and mediator model could be if living in a remote location increases the probability of both living alone and having lowered consciousness upon arrival to hospital. If on the other hand people who live in urban areas are more likely to be living alone and their proximity to the hospital decreases the probability of having lowered consciousness upon arrival to hospital then omitting this factor would induce a negative correlation between the error terms in the exposure and mediator model. The  $\widehat{NIE}(\boldsymbol{x})$ appears to be less sensitive to mediator-outcome confounding, as can be seen from Figure \ref{fig:sensanMean}b, where $0 \notin 95\%$ UI for $\rho_{\eta\xi}\in\left(-1,0.64 \right)$, and the effect is reversed, i.e. $NIE(\boldsymbol{x})_{u,\rho_{\eta\xi}}<0$, if $\rho_{\eta\xi}\in\left(0.74,1 \right)$. Finally, the $NIE(\boldsymbol{x})$ is not sensitive to exposure-outcome confounding (Figure \ref{fig:sensanMean}c). 
				
Figure \ref{fig:summary} summarizes the results of the sensitivity analyses for all the significant effects in Table \ref{table:effects} through ranges of $\rho_{\varepsilon\eta}$ (unobserved exposure-mediator confounding), $\rho_{\eta\xi}$ (mediator-outcome confounding) and $\rho_{\varepsilon\xi}$ (exposure-outcome confounding) that would render the effects significant and in the same direction as in Table \ref{table:effects} (light gray), non-significant (medium gray) or reversed (dark gray) at the 5\% level. These ranges coincide for the $NIE(\boldsymbol{x})$ for all three covariate patterns (no university, ages 60.4, 72.2 and 84.1 second panel from the top in Figure \ref{fig:summary}), i.e. the ranges seen from Figure \ref{fig:sensanMean}. The marginal natural direct effect (top-most panel) had similar results, with the effect rendered non-significant for $\rho_{\varepsilon\eta}\in\left(0.02,0.12 \right)$ or $\rho_{\eta\xi}\in\left(0.63,0.75 \right)$ and reversed for $\rho_{\varepsilon\eta}\in\left(0.12,1 \right)$ or $\rho_{\eta\xi}\in\left(0.75,1 \right)$.
				
For the $NDE(\boldsymbol{x})$ ranges coincide for ages 60.4 and 72.2 within a given level of education (no university, third panel from the top, and university, fourth panel from the top, Figure \ref{fig:summary}) and differ only slightly between education levels. For a patient without university education the effect is no longer significant for $\rho_{\eta\xi}\in\left(-1,-0.9 \right)$, the corresponding range for a patient with university education is $\rho_{\eta\xi}\in\left(-1,-0.88 \right)$, i.e. it takes unobserved mediator-outcome confounding that induces a strong negative correlation to render the natural direct effect non-significant. The natural direct effect is more sensitive to exposure-outcome confounding where the effect is no longer significant for $\rho_{\varepsilon\xi}\in\left(0.02,0.14 \right)$ and reversed for $\rho_{\varepsilon\xi}\in\left(0.14,1 \right)$ for both patients with and without university education. As previously stated unobserved exposure-outcome confounding could e.g. be due to pre-exposure socioeconomic factors not captured by education level. Low socioeconomic status has been linked to an increased risk of adverse outcome after stroke \cite{lindmark2014,cox2006,addo2012}. Depending on the effect of omitted socioeconomic factors on the probability of living alone, $\rho_{\varepsilon\xi}$ may be either positive (if the omitted factors increase the probability of living alone) or negative (the omitted factors decrease the probability of living alone). 
				
For an 84.1 year old patient the results are again quite similar with and without university education (two bottom-most panels of Figure \ref{fig:summary}) where the effect ceases to be significant for $\rho_{\eta\xi}\in\left(0.41,1 \right)$, for a patient without university education and $\rho_{\eta\xi}\in\left(0.42,1 \right)$ for a patient with university education. A plausible scenario here is that an unobserved genetic factor would have the same effect on the probability of having lowered consciousness upon arrival and of being dead or dependent at 3 months, either increasing or decreasing both, and thus  $\rho_{\eta\xi}$ is more likely to be positive than negative. Again, the direct effect appears to be more sensitive to exposure-outcome confounding than to mediator-outcome confounding, the effect is no longer significant for $\rho_{\varepsilon\xi}\in\left(-0.18,-0.02 \right)$ and reversed for $\rho_{\varepsilon\xi}\in\left(-1,-0.18 \right)$ for both patients with and without university education.

\section{Conclusion}
\label{sec:discussion}
The estimation of direct and indirect effects relies on strong assumptions about unconfoundedness. These assumptions are not testable using the observed data and so it is crucial that a mediation analysis be accompanied by a sensitivity analysis of the resulting estimates. We propose a sensitivity analysis method for mediation analysis based on probit regression models for both the mediator and the outcome. The sensitivity parameters introduced consist of the correlation between the error terms of the mediator and outcome models, as well as the correlation between the error terms of the mediator model and the model for the exposure assignment mechanism and the correlation between the error terms in the outcome model and the exposure assignment model. Incorporating these correlations into the estimation of the regression parameters allows us to obtain e.g. identification sets for the natural direct and indirect effects for a range of plausible correlation values. Sampling variability can be taken into account through the construction of uncertainty intervals. Our approach is able to take into account not only the mediator-outcome confounding that has been the focus of previous approaches but also exposure-mediator and exposure-outcome confounding. In addition, our method covers the situation where both the mediator and outcome are binary.

Using data from Riksstroke we performed a sensitivity analysis of the results from a mediation analysis of the effect of living alone on the probability of death or being dependent in ADL 3 months after stroke, with stroke severity (level of consciousness) upon arrival to hospital as mediator. In this study we did not have access to a rich set of pre-exposure covariates to adjust for and thus it was essential to perform a sensitivity analysis for not only mediator-outcome but also exposure-mediator and exposure-outcome confounding. The results of the sensitivity analysis were that the natural indirect effect was more sensitive to unobserved exposure-mediator confounding than to unobserved mediator-outcome confounding and not sensitive to unobserved exposure-outcome confounding. The natural direct effect was quite sensitive to unobserved exposure-outcome confounding, less sensitive to unobserved mediator-outcome confounding and not sensitive to unobserved exposure-mediator confounding. 
	
Although the method presented here is based on binary probit regression models for the exposure, mediator and outcome, this approach can be adapted to continuous exposures, mediators and/or outcomes. The method evaluates sensitivity to unobserved exposure-mediator, mediator-outcome, and exposure-outcome confounding separately. It would also be of interest to extend the method to investigation of the simultaneous effect of several types of unobserved confounding. Since a drawback to the method is its reliance on specifying parametric models for the exposure assignment mechanism, mediator and outcome, future work should also include generalizing it to semi-parametric mediation analysis which is less sensitive to model misspecification \cite{tchetgen2012,huber2014}. 

\section*{Appendix A: Probit based expressions for the total direct effect and pure indirect effect}

The total direct effect is defined as
$$NDE^* =\mathbb{E}\left[Y_i\left(1,M_i(1)\right)-Y_i\left(0,M_i(1)\right)\right],$$
and the pure indirect effect as
$$NIE^* =\mathbb{E}\left[Y_i\left(0,M_i(1)\right)-Y_i\left(0,M_i(0)\right)\right].$$
	
Assuming models \eqref{mstar} and \eqref{ystar} for the mediator and outcome and substituting \eqref{eq:expMprobit} and \eqref{eq:expYprobit} into the equivalent identification results to \eqref{eq:thm1nde}-\eqref{eq:thm1nie} outlined at the end of Section \ref{sec:identification} yields the following expressions for the conditional effects 
	\begin{align}
	NDE^*(\boldsymbol{x}) =&\left\{\Phi\left(\theta_0+\theta_1+\left( \boldsymbol{\theta}_4^\top+\boldsymbol{\theta}_5^\top\right) \boldsymbol{x}\right) - \Phi\left(\theta_0+\boldsymbol{\theta}_4^\top\boldsymbol{x}\right)\right\}\left(1-\Phi\left(\beta_0+\beta_1+\left( \boldsymbol{\beta}_2^\top+\boldsymbol{\beta}_3^\top\right) \boldsymbol{x}\right)\right) + \notag\\
	&\left\{\Phi\left(\theta_0+\theta_1+\theta_2+\theta_3+\left(\boldsymbol{\theta}_4^\top+\boldsymbol{\theta}_5^\top+\boldsymbol{\theta}_6^\top+\boldsymbol{\theta}_7^\top\right)\boldsymbol{x}\right) -\Phi\left(\theta_0+\theta_2+\left(\boldsymbol{\theta}_4^\top+\boldsymbol{\theta}_6^\top\right)\boldsymbol{x}\right)\right\}\times\notag\\ &\;\;\;\Phi\left(\beta_0+\beta_1+\left( \boldsymbol{\beta}_2^\top+\boldsymbol{\beta}_3^\top\right) \boldsymbol{x}\right), \notag \\[7pt] 
	NIE^*(\boldsymbol{x}) =&\left\{\Phi\left(\theta_0+\theta_2+\left(\boldsymbol{\theta}_4^\top+\boldsymbol{\theta}_6^\top\right)\boldsymbol{x}\right) - \Phi\left(\theta_0+\boldsymbol{\theta}_4^\top\boldsymbol{x}\right)\right\}\left\{\Phi\left(\beta_0+\beta_1+\left(\boldsymbol{\beta}_2^\top+\boldsymbol{\beta}_3^\top\right)\boldsymbol{x}\right)- \Phi\left(\beta_0+\boldsymbol{\beta}_2^\top\boldsymbol{x}\right)\right\}, \notag
	\end{align}
	and the $NDE^*$ and $NIE^*$ can be obtained by averaging the conditional effects over the population.

\section*{Appendix B: Standard errors of the estimators of the natural direct and indirect effects}	
The estimators $\widehat{NDE}(\boldsymbol{x})$ and $\widehat{NDE}(\boldsymbol{x})$ from \eqref{NDEhat} and \eqref{NIEhat}) are functions of the ML estimators  $\hat{\boldsymbol{\beta}}$ (or $\hat{\boldsymbol{\beta}}(\tilde{\rho}_{\varepsilon\eta})$, $\hat{\boldsymbol{\beta}}(\tilde{\rho}_{\eta\xi})$) of $\boldsymbol{\beta}=\left( \beta_0,\beta_1,\boldsymbol{\beta}_2,\boldsymbol{\beta}_3\right) $ and $\hat{\boldsymbol{\theta}}$ (or $\hat{\boldsymbol{\theta}}(\tilde{\rho}_{\eta\xi})$, $\hat{\boldsymbol{\theta}}(\tilde{\rho}_{\varepsilon\xi})$) of $\boldsymbol{\theta}=\left( \theta_0,\theta_1,\theta_2,\theta_3,\boldsymbol{\theta}_4,\boldsymbol{\theta}_5,\boldsymbol{\theta}_6,\boldsymbol{\theta}_7\right)$. From the delta method we have that $$\sqrt{n}\left( \widehat{NDE}(\boldsymbol{x})-NDE(\boldsymbol{x})\right) \xrightarrow{d} N\left( 0,\boldsymbol{\Lambda}\boldsymbol{\Sigma}\boldsymbol{\Lambda}^\top\right),$$
where $$\boldsymbol{\Sigma}=\left[ \begin{array}{cc}
\boldsymbol{\Sigma}_{\hat{\beta}} & \mathbf{0}\\
\mathbf{0} & \boldsymbol{\Sigma}_{\hat{\theta}}
\end{array}\right] ,$$ and $\boldsymbol{\Sigma}_{\hat{\beta}}$ and $\boldsymbol{\Sigma}_{\hat{\theta}}$ are the covariance matrices of  $\hat{\boldsymbol{\beta}}$ ($\hat{\boldsymbol{\beta}}(\tilde{\rho}_{\varepsilon\eta})$) and $\hat{\boldsymbol{\theta}}$ ($\hat{\boldsymbol{\theta}}(\tilde{\rho}_{\eta\xi})$, $\hat{\boldsymbol{\theta}}(\tilde{\rho}_{\varepsilon\xi})$), respectively, obtained through the inverse of the Fisher information matrices. The standard error of $\widehat{NDE}(\boldsymbol{x})$ is given by $\sqrt{\boldsymbol{\Lambda}\boldsymbol{\Sigma}\boldsymbol{\Lambda}^\top}$, where  $\boldsymbol{\Lambda}=\left( d_1,d_2,\boldsymbol{d}_3,\boldsymbol{d}_4,d_5,d_6,d_7,d_8,\boldsymbol{d}_9,\boldsymbol{d}_{10},\boldsymbol{d}_{11},\boldsymbol{d}_{12}\right) $ is the vector of partial derivatives of $NDE(\boldsymbol{x})$ wrt $\boldsymbol{\beta}$ and $\boldsymbol{\theta}$. Let
\begin{align*}
A&=\left\{\Phi\left(\theta_0+\theta_1+\left( \boldsymbol{\theta}_4^\top+\boldsymbol{\theta}_5^\top\right) \boldsymbol{x}\right) - \Phi\left(\theta_0+\boldsymbol{\theta}_4^\top\boldsymbol{x}\right)\right\} \\
B&=\left\{\Phi\left(\theta_0+\theta_1+\theta_2+\theta_3+\left(\boldsymbol{\theta}_4^\top+\boldsymbol{\theta}_5^\top+\boldsymbol{\theta}_6^\top+\boldsymbol{\theta}_7^\top\right)\boldsymbol{x}\right) -\Phi\left(\theta_0+\theta_2+\left(\boldsymbol{\theta}_4^\top+\boldsymbol{\theta}_6^\top\right)\boldsymbol{x}\right)\right\}\\
C&=\Phi\left( \beta_0+\boldsymbol{\beta}^\top_2\boldsymbol{x}\right)\\
D&=\phi\left(\theta_0+\theta_1+\theta_2+\theta_3+\left(\boldsymbol{\theta}_4^\top+\boldsymbol{\theta}_5^\top+\boldsymbol{\theta}_6^\top+\boldsymbol{\theta}_7^\top\right)\boldsymbol{x}\right),
\end{align*}
which gives
\begin{align*}
&d_1=A\left\lbrace  -\phi\left( \beta_0+\boldsymbol{\beta}^\top_2\boldsymbol{x}\right)\right\rbrace+B\phi\left( \beta_0+\boldsymbol{\beta}^\top_2\boldsymbol{x}\right)   \\
&d_2=0\\
&\boldsymbol{d}_3=d_1\boldsymbol{x}\\
&\boldsymbol{d}_4=\mathbf{0}\\
&d_5=\left\{\phi\left(\theta_0+\theta_1+\left( \boldsymbol{\theta}_4^\top+\boldsymbol{\theta}_5^\top\right) \boldsymbol{x}\right) - \phi\left(\theta_0+\boldsymbol{\theta}_4^\top\boldsymbol{x}\right)\right\}\left( 1-C\right) + \left\{D-\phi\left(\theta_0+\theta_2+\left(\boldsymbol{\theta}_4^\top+\boldsymbol{\theta}_6^\top\right)\boldsymbol{x}\right)\right\} C\\
&d_6=\phi\left(\theta_0+\theta_1+\left( \boldsymbol{\theta}_4^\top+\boldsymbol{\theta}_5^\top\right) \boldsymbol{x}\right)\left( 1-C\right)+DC\\
&d_7=\left\{D-\phi\left(\theta_0+\theta_2+\left(\boldsymbol{\theta}_4^\top+\boldsymbol{\theta}_6^\top\right)\boldsymbol{x}\right)\right\}C\\
&d_8=DC\\
&\boldsymbol{d}_9=d_5\boldsymbol{x}\\
&\boldsymbol{d}_{10}=d_6\boldsymbol{x}\\
&\boldsymbol{d}_{11}=d_7\boldsymbol{x}\\
&\boldsymbol{d}_{12}=d_8\boldsymbol{x}.
\end{align*}

For the natural indirect effect we have
$$\sqrt{n}\left( \widehat{NIE}(\boldsymbol{x})-NIE(\boldsymbol{x})\right) \xrightarrow{d} N\left( 0,\boldsymbol{\Gamma}\boldsymbol{\Sigma}\boldsymbol{\Gamma}^\top\right),$$
with $\boldsymbol{\Sigma}$ as before. The standard error of $\widehat{NIE}(\boldsymbol{x})$ is given by $\sqrt{\boldsymbol{\Gamma}\boldsymbol{\Sigma}\boldsymbol{\Gamma}^\top}$, where $\boldsymbol{\Gamma}= \left( g_1,g_2,\boldsymbol{g}_3,\boldsymbol{g}_4,g_5,g_6,g_7,g_8,\boldsymbol{g}_9,\boldsymbol{g}_{10},\boldsymbol{g}_{11},\boldsymbol{g}_{12}\right)$ the vector of partial derivatives of $NIE(\boldsymbol{x})$ wrt $\boldsymbol{\beta}$ and $\boldsymbol{\theta}$. Let
\begin{align*}
F&=\left\{\Phi\left(\theta_0+\theta_1+\theta_2+\theta_3+\left(\boldsymbol{\theta}_4^\top+\boldsymbol{\theta}_5^\top+\boldsymbol{\theta}_6^\top+\boldsymbol{\theta}_7^\top\right)\boldsymbol{x}\right) - \Phi\left(\theta_0+\theta_1+\left( \boldsymbol{\theta}_4^\top+\boldsymbol{\theta}_5^\top\right) \boldsymbol{x}\right)\right\} \\
G&=\left\{\Phi\left(\beta_0+\beta_1+\left(\boldsymbol{\beta}_2^\top+\boldsymbol{\beta}_3^\top\right)\boldsymbol{x}\right) -\Phi\left(\beta_0+\boldsymbol{\beta}_2^\top\boldsymbol{x}\right)\right\}, 
\end{align*}
which gives
\begin{align*}
&g_1=F\left\{\phi\left(\beta_0+\beta_1+\left(\boldsymbol{\beta}_2^\top+\boldsymbol{\beta}_3^\top\right)\boldsymbol{x}\right) -\phi\left(\beta_0+\boldsymbol{\beta}_2^\top\boldsymbol{x}\right)\right\} \\
&g_2=F\phi\left(\beta_0+\beta_1+\left(\boldsymbol{\beta}_2^\top+\boldsymbol{\beta}_3^\top\right)\boldsymbol{x}\right)\\
&\boldsymbol{g}_3=g_1\boldsymbol{x}\\
&\boldsymbol{g}_4=g_2\boldsymbol{x}\\
&g_5=\left\{D-\phi\left(\theta_0+\theta_1+\left( \boldsymbol{\theta}_4^\top+\boldsymbol{\theta}_5^\top\right)\boldsymbol{x}\right)\right\}G\\
&g_6=g_5\\
&g_7=DG\\
&g_8=g_7\\
&\boldsymbol{g}_9=g_5\boldsymbol{x}\\
&\boldsymbol{g}_{10}=\boldsymbol{g}_9\\
&\boldsymbol{g}_{11}=g_7\boldsymbol{x}\\
&\boldsymbol{g}_{12}=\boldsymbol{g}_{11}.
\end{align*}

	For the marginal effects we have that
	
	$$\sqrt{n}\left( \widehat{NDE}-NDE\right) \xrightarrow{d} N\left( 0,\boldsymbol{H}\boldsymbol{\Sigma}\boldsymbol{H}^\top\right),$$
	and
	$$\sqrt{n}\left( \widehat{NIE}-NIE\right) \xrightarrow{d} N\left( 0,\boldsymbol{K}\boldsymbol{\Sigma}\boldsymbol{K}^\top\right),$$
	where  $\widehat{NDE}=\frac{1}{n}\sum_{i=1}^n\widehat{NDE}(\boldsymbol{x}_i)$ and $\widehat{NIE}=\frac{1}{n}\sum_{i=1}^n\widehat{NIE}(\boldsymbol{x}_i)$. 
	
	The standard error of $\widehat{NDE}$ is given by $\sqrt{\boldsymbol{H}\boldsymbol{\Sigma}\boldsymbol{H}^\top}$, where  $\boldsymbol{H}=(h_1,h_2,\boldsymbol{h}_3,\boldsymbol{h}_4,h_5,h_6,h_7,h_8,\boldsymbol{h}_9,\boldsymbol{h}_{10},\boldsymbol{h}_{11},\boldsymbol{h}_{12})$ is the vector of partial derivatives of $NDE=\frac{1}{n}\sum_{i=1}^nNDE(\boldsymbol{x}_i)$ wrt $\boldsymbol{\beta}$ and $\boldsymbol{\theta}$, obtained by averaging the corresponding elements of $\boldsymbol{\Lambda}$.  
	
	The standard error of $\widehat{NIE}$ is given by $\sqrt{\boldsymbol{K}\boldsymbol{\Sigma}\boldsymbol{K}^\top}$, where  $\boldsymbol{K}=(k_1,k_2,\boldsymbol{k}_3,\boldsymbol{k}_4,k_5,k_6,k_7,k_8,\boldsymbol{k}_9,\boldsymbol{k}_{10},\boldsymbol{k}_{11},\boldsymbol{k}_{12})$ is the vector of partial derivatives of $NIE=\frac{1}{n}\sum_{i=1}^nNIE(\boldsymbol{x}_i)$ wrt $\boldsymbol{\beta}$ and $\boldsymbol{\theta}$, obtained by averaging the corresponding elements of $\boldsymbol{\Gamma}$.

\begin{acks}
The study was supported by the Swedish research council (grant no: 2012-5934). We are grateful to Ingeborg Waernbaum and Minna Genb\"{a}ck for comments that improved the paper, and to Riksstroke and the participating hospitals.
\end{acks}

\end{document}